\documentclass{article}
\usepackage[margin=40pt]{geometry}
\usepackage[T1]{fontenc}  
\usepackage{cite}
\usepackage{graphicx}
\usepackage{textcomp}
\usepackage{amsmath,amssymb,amsfonts,amsthm}
\usepackage{algorithmic}
\usepackage[ruled,linesnumbered]{algorithm2e}
\usepackage{textcomp}
\usepackage{hyperref}       % hyperlinks
\usepackage{url}            % simple URL typesetting
\usepackage{booktabs}       % professional-quality tables
\usepackage{nicefrac}
\usepackage{microtype}      % microtypography
\usepackage{lipsum}
\usepackage{enumitem}
\usepackage{listings}
\usepackage{graphicx}
\usepackage{epstopdf}
\usepackage{array}
\usepackage{subcaption}
\usepackage{wrapfig}
\usepackage{mathtools,cuted}
\usepackage{color}
\usepackage{colortbl}
\usepackage{cases}

\usepackage{stackengine,graphicx}
\stackMath

\usepackage{amssymb}
\usepackage[normalem]{ulem}
\usepackage{bm}
\usepackage{mathtools}
\usepackage{mathrsfs}
\usepackage{verbatim}

\usepackage{tikz}
\usetikzlibrary{shapes,arrows,chains}
\usepackage{authblk}

\newcommand{\keywords}[1]{%
  \vspace{1em}
  \noindent\textbf{Keywords: }#1
  \vspace{1em}
}

\begin{document}
\title{Task-Adaptive Reconstruction for Undersampled Multi-Coil MRI via Optimization-Guided Meta Learning}
\author[1]{Merham Fouladvand}
\author[2]{Peuroly Batra}
\affil[1,2]{Lincoln University,  Lincoln, Canterbury, New Zealand}
\maketitle

\begin{abstract}
    We propose a unified deep meta-learning framework for accelerated magnetic resonance imaging (MRI) that jointly addresses multi-coil reconstruction and cross-modality synthesis. Motivated by the limitations of conventional methods in handling undersampled data and missing modalities, we introduce a task-adaptive optimization network where
meta-learned parameters control the reconstruction dynamics
across varying acquisition conditions. Each phase of the network mimics a step of an adaptive forward-backward scheme with extrapolation, enabling the model to incorporate both data fidelity and nonconvex regularization in a principled manner. To enhance generalization across different acquisition settings, we integrate meta-learning, which enables the model to rapidly adapt to unseen sampling patterns and modality combinations using task-specific meta-knowledge. The proposed method is evaluated on the open source datasets, showing significant improvements in PSNR and SSIM over conventional supervised learning, especially under aggressive undersampling and domain shifts. Our results demonstrate the synergy of unrolled optimization, task-aware meta-learning, and modality fusion, offering a scalable and generalizable solution for real-world clinical MRI reconstruction.

\end{abstract}


\section{Introduction}
MRI reconstruction from undersampled measurements remains a
challenging inverse problem. While parallel imaging and compressed
sensing have significantly reduced acquisition time, aggressive
undersampling still leads to aliasing artifacts and unstable
reconstruction. Recent deep learning approaches have shown strong
performance but often require retraining when acquisition conditions
change. Despite these advantages, MRI examinations often require long acquisition times, which may lead to patient discomfort, motion artifacts, and limited scanner throughput. To mitigate these limitations, accelerated acquisition techniques have been developed to reduce the amount of measured data while maintaining acceptable reconstruction quality.

Parallel Imaging is widely used to accelerate imaging by employing multiple receiver coils to acquire signals simultaneously. Each coil observes the underlying image through a spatial sensitivity profile, allowing partial sampling of k-space while preserving sufficient information for reconstruction. However, recovering a high-quality image from undersampled measurements remains an ill-posed inverse problem. The reconstruction process must balance fidelity to the acquired data with prior knowledge about the structure of medical images.

In parallel with advances in reconstruction techniques, clinical MRI increasingly relies on multi-modality imaging. Different MRI contrasts, such as T1-weighted, T2-weighted, and FLAIR images, provide complementary information about tissue composition and pathology. Combining information from multiple modalities can significantly improve diagnostic interpretation. However, acquiring all modalities for every patient is often impractical due to time constraints, patient motion, and acquisition artifacts. Consequently, there is growing interest in methods that can jointly reconstruct undersampled data and infer missing modalities.

Deep learning has recently shown strong potential for solving inverse problems in medical imaging. In particular, algorithm-unrolled neural networks have emerged as a promising framework for integrating domain knowledge with data-driven learning. These models interpret iterative optimization algorithms as neural network architectures, where each iteration corresponds to a learnable computational block. This structure allows the model to maintain physical consistency with the imaging process while learning effective image priors from data.

Although unrolled networks have achieved encouraging results in MRI reconstruction, two limitations remain. First, many existing approaches focus solely on single-modality reconstruction and do not explicitly exploit relationships between different imaging contrasts. Second, models trained for a specific sampling pattern or acquisition protocol often struggle to generalize to new scanning conditions.

To address these challenges, this work introduces a unified reconstruction framework that combines algorithm unrolling with meta-learning. The proposed approach simultaneously handles multi-coil reconstruction and cross-modality synthesis from undersampled measurements. By incorporating meta-learning into the reconstruction process, the model learns transferable knowledge that enables rapid adaptation to new acquisition patterns and modality configurations. In addition, the underlying optimization algorithm provides theoretical convergence guarantees, ensuring stability of the iterative reconstruction process.

The main contributions of this work can be summarized as follows:\\

1. Unified Reconstruction and Modality Synthesis.

We propose a framework that jointly reconstructs multi-coil MR images and synthesizes missing modalities within a single optimization-driven architecture.\\

2. Algorithm-Unrolled Network with Adaptive Optimization. 

The proposed model is derived from an adaptive forward–backward optimization scheme with extrapolation, enabling each network phase to correspond to a principled optimization step.\\

3. Meta-Learning for Cross-Task Adaptation.

A meta-learning mechanism is introduced to improve generalization across varying sampling patterns, acceleration factors, and modality combinations.\\

4. Convergence-Aware Network Design.

The reconstruction process is guided by a provably convergent optimization algorithm, providing theoretical support for the stability of the network.\\

Experimental results on publicly available datasets demonstrate that the proposed framework achieves improved reconstruction accuracy and robustness compared with conventional supervised learning approaches, particularly under aggressive undersampling and domain shifts.

\section{Related Work}\label{related}
Magnetic Resonance Imaging (MRI) reconstruction has been extensively studied through both traditional optimization-based methods and modern learning-based techniques. Classical parallel MRI (pMRI) methods such as SENSE\cite{pruessmann1999sense} and GRAPPA\cite{griswold2002generalized} utilize multiple receiver coils to accelerate image acquisition, relying heavily on accurate coil sensitivity estimation. However, their performance can deteriorate under high undersampling rates or in the presence of noise.

Recent works have explored deep learning methods to improve reconstruction quality\cite{bian2020deep}. For instance, deep cascaded networks and variational networks unroll iterative optimization algorithms into trainable neural architectures, enabling data-driven learning of regularization priors. Several studies have also incorporated domain knowledge into the network design, such as incorporating the Fourier transform to enforce k-space consistency or using complex-valued networks to handle MRI data more naturally.

In the context of multi-modality MRI\cite{usman2017brain,bian2022learnable}, prior research has primarily focused on synthesizing missing modalities from fully sampled data using conditional GANs\cite{wang2018high,quan2018compressed} or paired encoder-decoder architectures. While these approaches generate plausible images, they typically require fully scanned source modalities, limiting their practical utility in accelerated MRI acquisition. A variational framework\cite{bian2022learnable} that jointly reconstructs and synthesizes MR images across multiple modalities. Their model incorporates modality-specific encoders, a shared decoder, and learned variational priors, enabling efficient synthesis from undersampled inputs. This approach inspired the joint modeling adopted in our current work, where both fusion and reconstruction are handled within a unified deep learning structure. For instance, Hu et al.~\cite{hu2020supervised} proposed a supervised descent learning method for directional electromagnetic logging-while-drilling (LWD) inversion, showing that a learning-based approach can effectively solve complex inverse problems in geophysical settings. Complementing this, a theory-guided deep neural network framework was introduced for time-domain electromagnetic (TDEM) inversion and simulation using differentiable programming~\cite{hu2021theory}, offering a robust hybridization of physical modeling and data-driven inference.

Deep learning-based MRI reconstruction methods have rapidly advanced in recent years, particularly in addressing the challenges of undersampled data, multi-coil imaging, and multi-modality synthesis. Traditional approaches to MRI reconstruction often rely on fixed priors and hand-crafted regularizations. In contrast, recent methods focus on integrating data-driven models with domain knowledge, yielding improved flexibility and performance.

In the domain of multi-coil pMRI\cite{bian2022optimal,zhou2019multi}, several methods proposed an optimal control\cite{WCP18,BC17} formulation is very useful. By jointly reconstructing channel-wise images using an end-to-end learnable model, their method addresses one of the core limitations in traditional pMRI workflows.

Self-supervised learning has also been applied to quantitative MRI. Several models\cite{shurrab2022self,usman2017brain,bian2024improving,wang2018high} demonstrated that self-supervision combined with physics-based model reinforcement can significantly improve the accuracy of quantitative mapping (e.g., T1 estimation) while reducing scan time. These insights support the utility of hybrid learning paradigms in constrained MRI settings.

More recently, a diffusion model for accelerated MRI\cite{HKK18} and quantitative MRI \cite{bian2024accelerating}reconstruction using domain-conditioned prior guidance. This method leverages generative modeling to learn sample-consistent priors conditioned on domain information, achieving state-of-the-art results in both structural fidelity and quantitative accuracy. 

Building upon these advancements, our proposed method adopts a deep unrolled network structure with learned fusion operators and provably convergent optimization, aiming to simultaneously address the challenges of undersampled k-space acquisition, modality synthesis, and multi-coil reconstruction.

The described approach advances these lines of work by proposing a deep unrolled architecture that jointly performs multi-coil reconstruction and cross-modality synthesis directly from undersampled data. We introduce a learnable fusion operator across modalities and propose a provably convergent algorithm with smoothed nonconvex regularization, bridging theoretical guarantees with empirical performance.

\section{Method}\label{method}

The goal of this work is to reconstruct high-quality MRI images and synthesize missing modalities from undersampled multi-coil k-space measurements. To achieve this objective, we design a learning framework derived from an optimization algorithm and implement it as a deep unrolled neural network.

Let $x$ denote the set of reconstructed images corresponding to different modalities, and let the observed measurements be obtained from multiple receiver coils under an undersampling pattern. Unlike conventional reconstruction models that employ a fixed
regularization term, our formulation allows the regularization
component to be dynamically adapted through meta-learned parameters.
This enables the reconstruction process to adjust its implicit prior
based on the acquisition task.

The detailed steps are described in Algorithm~\ref{alg:afbe}.

\begin{algorithm}
\caption{Adaptive Forward--Backward Reconstruction Scheme}
\label{alg:afbe}
\begin{algorithmic}[1]
\STATE \textbf{Input:} Initial estimate $x^{(0)}$, smoothing parameter $\eta^{(0)}>0$, extrapolation factor $\tau^{(0)}\in(0,1)$, decay constants $\sigma,\mu,\delta\in(0,1)$, stopping threshold $\epsilon>0$
\FOR{$k=0,1,\dots,K-1$}
    \STATE Compute a fidelity-driven update
    \[
    z^{(k)} = x^{(k)} - \alpha_k \nabla f\!\left(x^{(k)}\right).
    \]
    \STATE Compute a regularization-aware candidate
    \[
    u^{(k)} = z^{(k)} - \beta_k \nabla r_{\eta^{(k)}}\!\left(z^{(k)}\right).
    \]
    \IF{$\phi_{\eta^{(k)}}(u^{(k)}) \le \phi_{\eta^{(k)}}(x^{(k)}) - \frac{\delta}{2}\|u^{(k)}-x^{(k)}\|^2$}
        \STATE Set $\bar{x}^{(k+1)} = u^{(k)}$
    \ELSE
        \STATE Set
        \[
        \bar{x}^{(k+1)} = x^{(k)} - \gamma_k \nabla \phi_{\eta^{(k)}}\!\left(x^{(k)}\right).
        \]
    \ENDIF
    \STATE Form an extrapolated candidate
    \[
    y^{(k+1)} = \bar{x}^{(k+1)} + \tau^{(k)}\big(\bar{x}^{(k+1)}-x^{(k)}\big).
    \]
    \IF{$\phi_{\eta^{(k)}}(y^{(k+1)}) \le \phi_{\eta^{(k)}}(\bar{x}^{(k+1)})$}
        \STATE Accept $x^{(k+1)} = y^{(k+1)}$ and update
        \[
        \tau^{(k+1)} = \min\!\left\{\frac{\tau^{(k)}}{\mu},\,1\right\}
        \]
    \ELSE
        \STATE Accept $x^{(k+1)} = \bar{x}^{(k+1)}$ and update
        \[
        \tau^{(k+1)} = \mu \tau^{(k)}.
        \]
    \ENDIF
    \IF{$\|\nabla \phi_{\eta^{(k)}}(x^{(k+1)})\| < \sigma \eta^{(k)}$}
        \STATE Set $\eta^{(k+1)} = \sigma \eta^{(k)}$
    \ELSE
        \STATE Set $\eta^{(k+1)} = \eta^{(k)}$
    \ENDIF
    \IF{$\eta^{(k+1)} < \epsilon$}
        \STATE \textbf{break}
    \ENDIF
\ENDFOR
\STATE \textbf{Output:} $x^{(k+1)}$
\end{algorithmic}
\end{algorithm}

This adaptive algorithm provides a convergence-guaranteed framework and is especially well-suited for unrolling into deep network layers where each iteration corresponds to one network phase. The smoothed regularizer enables gradient-based training and handles nonconvexity in a controlled manner. The extrapolation and adaptive $\eta_k$ updates allow the model to navigate the optimization landscape efficiently, balancing convergence speed and stability.

The training objective is formulated as a composite loss that jointly
penalizes synthesis error, perceptual inconsistency, cross-domain
mapping discrepancy, and reconstruction error. The loss is defined as

\begin{subequations}\label{eq:training_loss}
\begin{align}
\mathcal{L}(x_{0},x_{1};\hat{x}_{0},\hat{x}_{1})
=&\; \|x_{0}-\hat{x}_{0}\|_{2}^{2}  \\
&+ \lambda_{1}\big(1-\mathrm{SSIM}(x_{0},\hat{x}_{0})\big)  \\
&+ \lambda_{2}\left\|
\tilde{\mathbf g}_{0}\!\left(
\mathcal{G}\big(\mathbf g_{1}(x_{1})\big)
\right)-\hat{x}_{0}
\right\|_{2}^{2}  \\
&+ \lambda_{3}\|x_{1}-\hat{x}_{1}\|_{2}^{2}.
\end{align}
\end{subequations}

Here, $x_{0}$ denotes the synthesized image produced by the synthesis
branch, while $x_{1}$ represents the reconstructed image obtained at the
$T$-th phase of the network. The corresponding reference images
are denoted by $\hat{x}_{0}$ and $\hat{x}_{1}$.

The loss incorporates four complementary terms. The first component
measures the squared reconstruction error between the synthesized image
and its reference. The second term introduces the structural similarity
index (SSIM) to preserve perceptual structure and contrast information.
The third component enforces a cycle-style consistency constraint by
mapping the reconstructed image through the learned transformation
$\tilde{\mathbf g}_{0}(\mathcal{G}(\mathbf g_{1}(\cdot)))$ and comparing
the resulting estimate with the synthesized target. The final term
penalizes the deviation between the reconstructed output and its
ground-truth image. The scalar weights $\lambda_{1},\lambda_{2},\lambda_{3}$
control the relative influence of the perceptual, consistency, and
reconstruction penalties during training.

\begin{figure}
\includegraphics[width=0.6\linewidth]{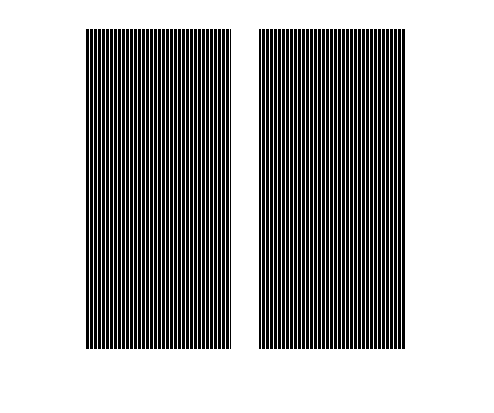}
\caption{Uniform Cartesian Mask. Sampling Ratio 31.56\% , this is the mask we used in the training. The paper we submitted also used this mask. I added more data and this is the newest result: Phase5, Avg REC PSNR is 41.7114 dB, Avg relative error is 0.026917 dB,  Avg ssim is 0.9719 dB}
\end{figure}  

\begin{figure}
\includegraphics[width=0.6\linewidth]{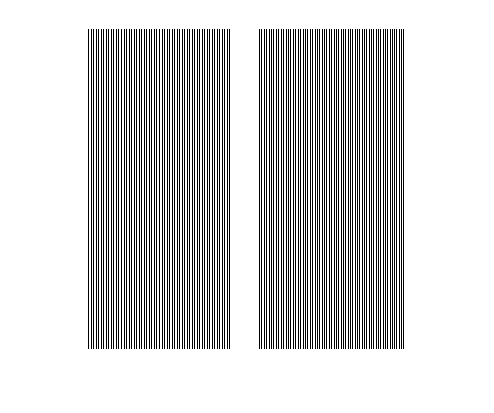}
\caption{Uniform Cartesian Mask. Sampling Ratio 62.81\% ,  Reconstruction uses the weights of phase 5 (same as fig1), Avg REC PSNR is 26.3405 dB, Avg relative error is 0.158550 dB,  Avg ssim is 0.8257}
\end{figure} 

\begin{figure}
\includegraphics[width=0.6\linewidth]{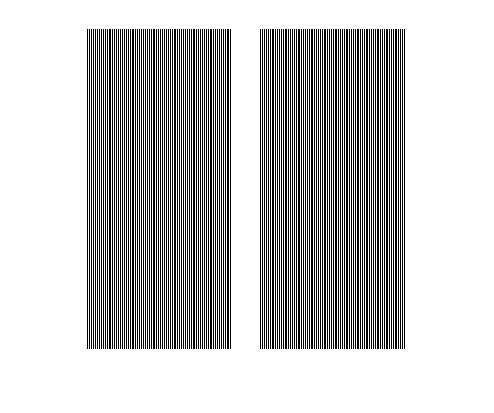}
\caption{Uniform Cartesian Mask. Sampling Ratio 51.56\% ,  Reconstruction uses the weights of phase 5 (same as fig1),  Avg REC PSNR is 29.8603 dB, Avg relative error is 0.105742 dB,  Avg ssim is 0.8982}
\end{figure} 

\begin{figure}
\includegraphics[width=0.6\linewidth]{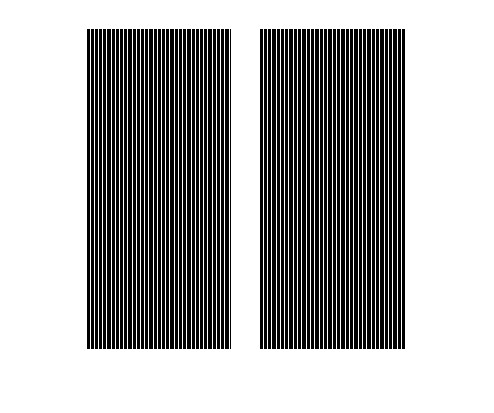}
\caption{Uniform Cartesian Mask. Sampling Ratio 30.63\% ,  Reconstruction uses the weights of phase 5 (same as fig1),  Avg REC PSNR is 34.2160 dB, Avg relative error is 0.063928 dB,  Avg ssim is 0.9083 }
\end{figure} 

\begin{figure}
\includegraphics[width=0.6\linewidth]{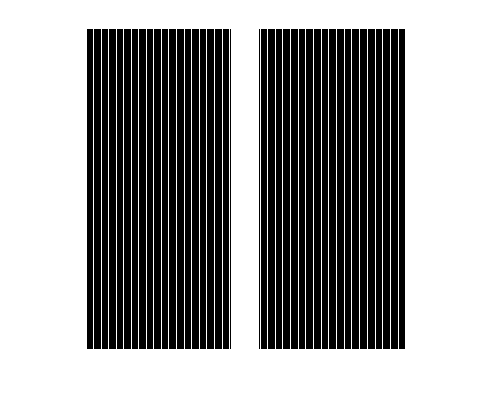}
\caption{Uniform Cartesian Mask. Sampling Ratio 21.25\% ,  Reconstruction uses the weights of phase 5 (same as fig1), Avg REC PSNR is 32.2787 dB, Avg relative error is 0.079973 dB,  Avg ssim is 0.8814 }

\end{figure} 

\begin{figure}
\includegraphics[width=0.6\linewidth]{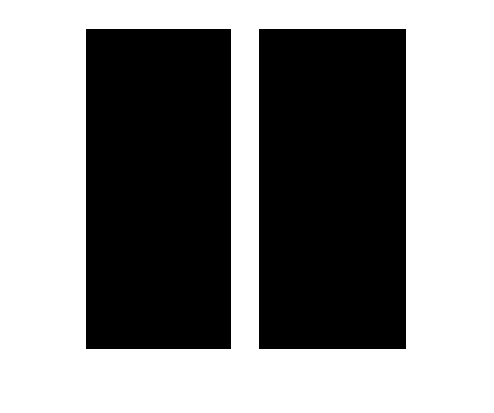}
\caption{Uniform Cartesian Mask. Sampling Ratio 8.75\% , Reconstruction uses the weights of phase 5 (same as fig1), Avg REC PSNR is 27.1134 dB, Avg relative error is 0.144922 dB,  Avg ssim is 0.7453 }

\end{figure} 

\begin{figure}
\includegraphics[width=0.6\linewidth]{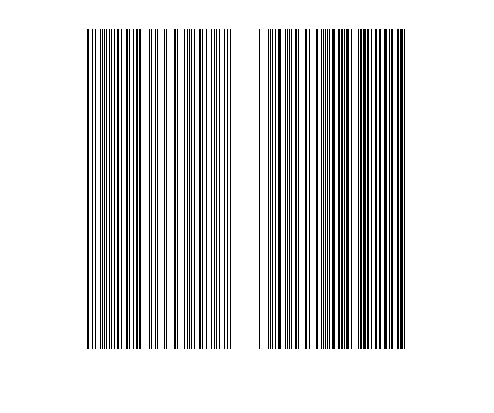}
\caption{Random Cartesian Mask. Sampling Ratio 50.31\% ,  Reconstruction uses the weights of phase 5 (same as fig1), Avg REC PSNR is 24.3860 dB, Avg relative error is 0.198565 dB,  Avg ssim is 0.7563 }
 
\end{figure}

\begin{figure}
\includegraphics[width=0.6\linewidth]{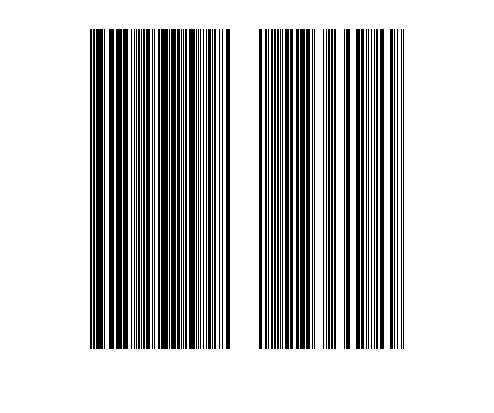}
\caption{Random Cartesian Mask. Sampling Ratio 50.31\% ,  Reconstruction uses the weights of phase 5 (same as fig1), Avg REC PSNR is 27.7152 dB, Avg relative error is 0.135180 dB,  Avg ssim is 0.8051 }

\end{figure}

\begin{figure}
\includegraphics[width=0.6\linewidth]{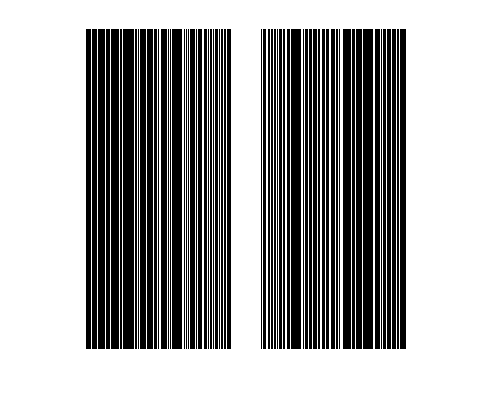}
\caption{Random Cartesian Mask. Sampling Ratio 31.88\% ,  Reconstruction uses the weights of phase 5 (same as fig1), Avg REC PSNR is 31.0304 dB, Avg relative error is 0.092377 dB,  Avg ssim is 0.8564 }

\end{figure}  

\begin{figure}
\includegraphics[width=0.6\linewidth]{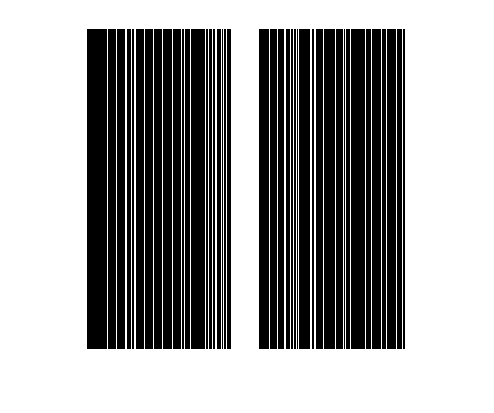}
\caption{Random Cartesian Mask. Sampling Ratio 23.44\% ,  Reconstruction uses the weights of phase 5 (same as fig1), Avg REC PSNR is 30.9118 dB, Avg relative error is 0.093482 dB,  Avg ssim is 0.8396 }

\end{figure} 

\begin{figure}
\includegraphics[width=0.6\linewidth]{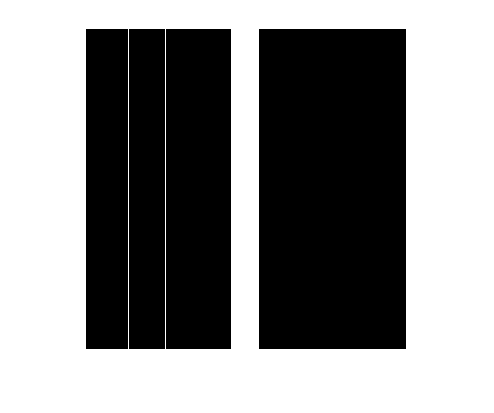}
\caption{Random Cartesian Mask. Sampling Ratio 9.38\% ,  Reconstruction uses the weights of phase 5 (same as fig1),Avg REC PSNR is 27.1566 dB, Avg relative error is 0.144207 dB,  Avg ssim is 0.7506 }

\end{figure} 

\begin{figure}
\includegraphics[width=0.6\linewidth]{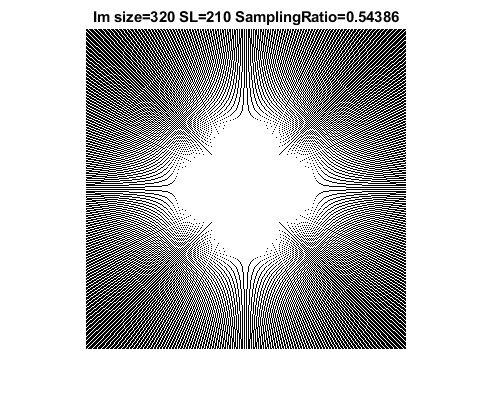}
\caption{Radio Mask. Sampling Ratio 54.386\% ,  Reconstruction uses the weights of phase 5 (same as fig1),  Avg REC PSNR is 21.3429 dB, Avg relative error is 0.282210 dB,  Avg ssim is 0.6598 }
\end{figure}

\begin{figure}
\includegraphics[width=0.6\linewidth]{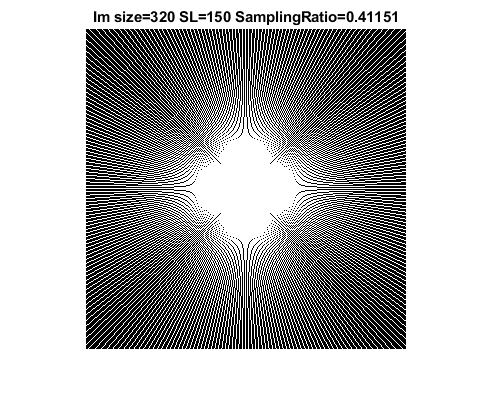}
\caption{Radio Mask. Sampling Ratio 41.15\% ,  Reconstruction uses the weights of phase 5 (same as fig1),  Avg REC PSNR is 21.6037 dB, Avg relative error is 0.273975 dB,  Avg ssim is 0.6530}
\label{fig3}
\end{figure}

\begin{figure}
\includegraphics[width=0.6\linewidth]{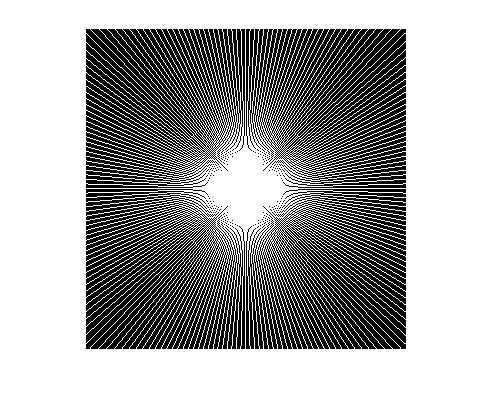}
\caption{Radio Mask. Sampling Ratio 31.288\% ,  Reconstruction uses the weights of phase 5 (same as fig1),  Avg REC PSNR is 21.9603 dB, Avg relative error is 0.262960 dB,  Avg ssim is 0.6464 }
\label{fig4}
\end{figure}  

\begin{figure}
\includegraphics[width=0.6\linewidth]{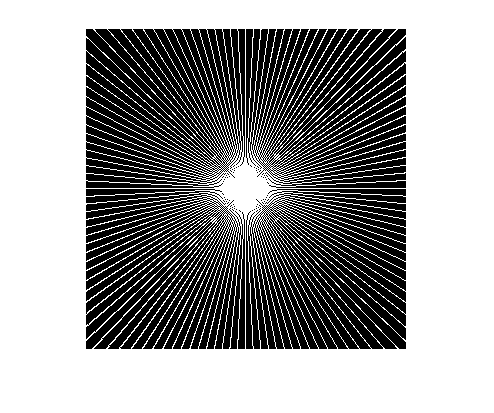}
\caption{Radio Mask. Sampling Ratio 21.40\% ,  Reconstruction uses the weights of phase 5 (same as fig1),  Avg REC PSNR is 23.2455 dB, Avg relative error is 0.226482 dB,  Avg ssim is 0.6487 }
\label{fig5}
\end{figure}

\begin{figure}
\includegraphics[width=0.6\linewidth]{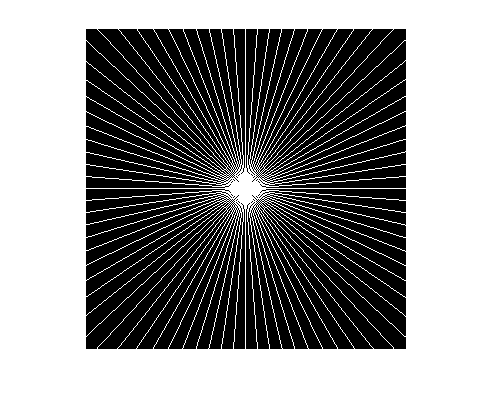}
\caption{Radio Mask. Sampling Ratio 12.67\% ,  Reconstruction uses the weights of phase 5 (same as fig1), Avg REC PSNR is 25.6671 dB, Avg relative error is 0.170824 dB,  Avg ssim is 0.6595}
\label{fig6}
\end{figure}

\begin{figure}
\includegraphics[width=0.6\linewidth]{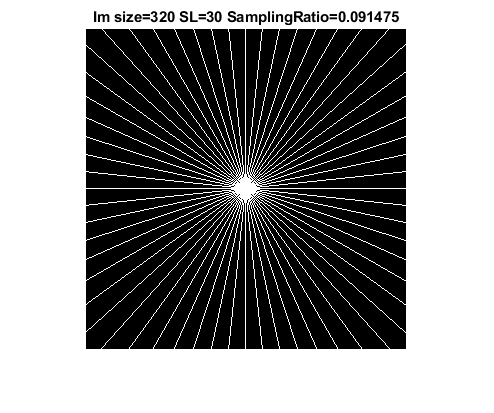}
\caption{Radio Mask. Sampling Ratio 9.14\% ,  Reconstruction uses the weights of phase 5 (same as fig1), Avg REC PSNR is 26.7459 dB, Avg relative error is 0.150878 dB,  Avg ssim is 0.6495}
\end{figure}

\begin{figure}
\includegraphics[width=0.6\linewidth]{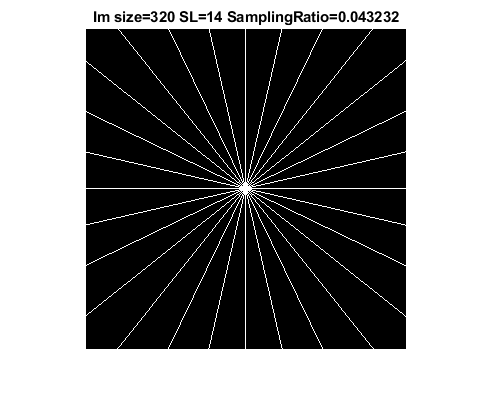}
\caption{Radio Mask. Sampling Ratio 4.32\% ,  Reconstruction uses the weights of phase 5 (same as fig1), Avg REC PSNR is 23.4165 dB, Avg relative error is 0.221283 dB,  Avg ssim is 0.5490}
\end{figure}

\section{Conclusion}
This study demonstrates how algorithm unrolling combined with adaptive optimization can bridge the gap between theoretical convergence guarantees and practical performance in data-driven MRI reconstruction. By drawing parallels with recent advancements in geophysical inverse problems and theory-guided neural networks, we highlight the broader applicability of our approach to a wide class of physics-informed inverse problems. The modularity of the proposed method also opens the door for future extensions, including integration with diffusion priors, domain-conditioned guidance, or meta-learned optimizers tailored to specific anatomy or scanner settings.

\bibliographystyle{IEEEtran}
\bibliography{ref}

\end{document}